\documentclass [11pt,fleqn] {article}
\usepackage{amsmath}
\usepackage{fleqn}
\usepackage{amssymb}
\DeclareMathAlphabet{\mathsl}{OT1}{cmr}{m}{sl}
\DeclareFontShape{OT1}{cmss}{m}{it}{ <-> ssub * cmss/m/sl }{}
\newcommand{\IN}{{\mathbb{N}}} 

\newcommand{\IF}{{\mathbb{F}}}
\newcommand{\IR}{{\mathbb{R}}}

\newcommand{\QCB}{\mathsf{QCB}}
\newcommand{\Seq}{\mathsf{Seq}}
\newcommand{\kHaus}{\mathsf{kHaus}}

\newcommand{\Top}{\mathsf{Top}}
\newcommand{\Equ}{\mathsf{Equ}}

\newcommand{\INNN}{\IN^{\IN^\IN}}

\newcommand{\id}{{\mathsl{id}}}
\newcommand{\w}[1]{{\mathtt{#1}}}
\newcommand{\RE}{{\mathsf{R}_{\mathrm{E}}}}
\newcommand{\RI}{{\mathsf{R}_{\mathrm{I}}}}

\newcommand{\Nk} [1]{\mathsf{N}\langle#1\rangle}

\newcommand{\REk}[1]{\mathsf{R}_{\mathrm{E}}\langle#1\rangle}
\newcommand{\RIk}[1]{\mathsf{R}_{\mathrm{I}}\langle#1\rangle}
\newcommand{\RSk}[1]{\mathsf{R}_{\mathrm{S}}\langle#1\rangle}

\newcommand{\WM}{{\mathsf{M}}}

\newcommand{\WMProd}{{\mathsf{M}_\mathrm{P}}}  

\newcommand{\WX}{{\mathsf{X}}}
\newcommand{\WY}{{\mathsf{Y}}}

\newcommand{\Wzwei}{{\mathsf{2}}}

\newcommand{\leins}{{\ell_1}}
\newcommand{\Tsum}{\mathop{\textstyle\sum}}

\newcommand{\Tprod}{\mathop{\textstyle\prod}}

\newcommand{\Tbigcap}{\mathop{\textstyle\bigcap}}

\newcommand{\nm}[1]{\|#1\|_1}

\newcommand{\setof}[1]{{\mid}#1{\mid}}

\newtheorem{Theorem}{Theorem}[section]
\newtheorem{Proposition}[Theorem]{Proposition}
\newtheorem{Lemma}[Theorem]{Lemma}
\newtheorem{Corollary}[Theorem]{Corollary}

\newenvironment{Proof}{\pagebreak[3]\noindent\emph{Proof.}}
{\par\hspace*{\fill}$\square$\pagebreak[3]\medskip}

{\end{list}}
\newenvironment{mycases}{\begin{list}{}{%
\setlength{\labelwidth}{2em}
\setlength{\labelsep}{5pt}
\setlength{\itemindent}{0em}%
\setlength{\leftmargin}{\labelwidth}%
\addtolength{\leftmargin}{\labelsep}%
\setlength{\rightmargin}{0em}%
\setlength{\topsep}{8pt plus 2pt minus 4pt}%
\setlength{\parskip}{0pt plus 1pt}%
\setlength{\itemsep}{4pt plus 2pt minus 1pt}%
\setlength{\parsep} {4pt plus 2pt minus 1pt}%
}}%
{\end{list}}
\begin{document}

\title{$\IN^{\IN^\IN}$ does not satisfy Normann's condition}

\author{Matthias Schr\"oder}
\date{\quad}

\maketitle

\begin{abstract}
We prove that the Kleene-Kreisel space $\IN^{\IN^\IN}$ does not satisfy
Normann's condition.
A topological space $\WX$ is said to fulfil Normann's condition,
if every functionally closed subset of $\WX$ is an intersection
of clopen sets.
The investigation of this property is motivated by its
strong relationship to a problem in Computable Analysis.
D.~Normann has proved that
in order to establish non-coincidence of
the extensional hierarchy and the intensional hierarchy
of functionals over the reals
it is enough to show that $\INNN$ fails the above condition.
\end{abstract}


\section{Introduction}

The Kleene-Kreisel continuous functionals over the natural numbers 
play an important role in mathematical logic
as well as in higher type computability
\cite{Kle:CountFunc,Krei:AnaFunc,Nor:LNM}.
A simple way of defining this hierarchy is to construct it
as a sequence of exponentials in an appropriate cartesian closed category
by applying the recursion formula $\Nk{0}:=\IN$ and $\Nk{k+1}:=\IN^{\Nk{k}}$.
In this paper we use the cartesian closed category $\QCB$ as our ambient category
\cite{Sim:towards}.
This full subcategory of $\Top$ has as objects all quotients
of countably based topological spaces.
Alternatives are the category $\Seq$ of sequential spaces
or the category $\kHaus$ of Hausdorff Kelley spaces~\cite{ELS:CCGS}.

The main goal of this paper is to prove that
the space $\INNN=\Nk{2}$ contains a functionally closed
subset\footnote{A subset $A$ of a topological space
  is \emph{functionally closed}, if $A$ is the preimage of $0$
  under a continuous function into the unit interval $[0;1]$.}
which can not be represented as an intersection of clopen sets.
Since $\INNN$ is hereditarily Lindel\"of, this is equivalent
to saying that the completely regular reflection of $\INNN$
is not zero-dimensional.
Note that $\INNN$ means the exponential formed in $\QCB$ 
(or equivalently in $\Seq$ or $\kHaus$, see \cite{ELS:CCGS})
to the basis $\IN$ and the exponent $\IN^\IN$.
So $\INNN$ is topologised by
the sequentialisation of the compact-open topology
on the set of continuous functions from the Baire space
to the discrete space $\IN$, see \cite{ELS:CCGS,Sch:phd}.
It is well-known that this topology is strictly finer than
the compact-open topology on $\INNN$.
From \cite{Sch:NNNnotT3} we know that $\INNN$ is neither zero-dimensional
nor regular.

As an important consequence of this result, we obtain
that the extensional hierarchy and the intensional hierarchy
of functionals over the real numbers do not coincide.
These two hierarchies have been introduced by 
Bauer, Escard{\'o} and Simpson
to model two approaches to higher type computation
over the real numbers in functional programming
\cite{BES:paradigms}.
Normann proved that the two hierarchies agree up to level $k+1$
if, and only if,
the Kleene-Kreisel space $\Nk{k}$ has the property that
functionally closed sets are intersections of clopen sets
(see \cite[Theorems 4.17 \& 5.5]{Nor:comparing}).
We therefore refer to this property as \emph{Normann's condition}.
Bauer, Escard{\'o} and Simpson had already observed 
that the coincidence question for level $3$
is related to topological properties of the space $\INNN$,
see \cite{BES:paradigms}.

In Section~\ref{sec:Definition:WM} we construct a Polish space $\WM$
that arises as the sequential coreflection of some zero-dimensional space,
but is not zero-dimensional itself.
So $\WM$ is a totally disconnected metric space that does not fulfil
Normann's condition.
In Section~\ref{sec:WM:into:INNN} we prove that $\WM$
is a retract of $\INNN$.
Both results combined entail that $\INNN$ does not
satisfy Normann's condition.
In Section~\ref{sec:application:CompAna} we briefly discuss
the extensional and the intensional hierarchies
of functionals over the reals.


\pagebreak[3]
\section{Definition of the Polish space $\WM$}
\label{sec:Definition:WM}

We will define $\WM$ to be a closed subspace
of the real vector space $\ell_1$.
The space $\ell_1$ consists of those elements $x \in \IR^\IN$
for which the $1$-norm $\nm{x}$ defined by
\[
  \nm{x} :=\sum_{i \in \IN} |x(i)|
\]
is less than $\infty$.
It is well-known that $(\ell_1,\nm{.})$ is a Banach space.
In abuse of notation, we henceforth denote by $\ell_1$
the countably based topological space
that carries the topology induced by the $\ell_1$-metric
$(x,y) \mapsto \nm{x -y}$.
We will later need the  following characterisation 
of convergence of sequences in $\ell_1$
which is folklore in functional analysis.

\pagebreak[3]
\begin{Lemma}\label{l:ConvInEll}
 A sequence $(x_n)_n$ converges in $\ell_1$ to $x_\infty$
 if, and only if, (a) and (b) hold:
 \begin{itemize}
  \item[(a)]
   For all $i \in \IN$, $(x_n(i))_n$ converges to $x_\infty(i)$ in $\IR$.
  \item[(b)] 
   The sequence $(\nm{x_n})_n$ converges to $\nm{x_\infty}$ in $\IR$.
 \end{itemize}
\end{Lemma}

\medskip

For the construction of $\WM$, we define for $i \in \IN$ the set $M_i$ by
\[
  M_i:=\big\{j \cdot 2^{-i} \,\big|\, j \in \{0,\dotsc,2^i\} \big\}
   =\{0,1\cdot 2^{-i},2\cdot 2^{-i},3\cdot 2^{-i}, \dotsc, (2^i-1) \cdot 2^{-i},1\}
  \,.    
\]  
Then the subspace $\WM$ of $\leins$ with underlying set
\begin{equation}
  |\WM|:= \big\{ x \in \Tprod_{i \in \IN} M_i : \nm{x} < \infty \big\}
\end{equation}
is closed in $\leins$.
This is due to the fact that the sets $M_i$ are closed in $\IR$.
So $\WM$ is a Polish space, with the restriction of the $\ell_1$-metric
being a complete metric for $\WM$.

We show in a similar way as in \cite[Example 6.2.19]{En89}
that $\WM$ is not zero-dimensional.

\begin{Lemma}\label{l:unitball:notZO}
 The unit ball
 $B_\WM(0^\omega;1):= \{ x \in \WM : \nm{x} <1 \}$
 does not contain any clopen neighbourhoud of
 the constant zero-function $0^\omega \in \leins$.
\end{Lemma}

\begin{Proof}
 Let $V$ be any open set with $0^\omega \in V \subseteq B_\WM(0^\omega;1)$.
 By recursion we construct a sequence $(a_k)_k \in \Tprod_{i \in \IN}M_i$
 such that
 \begin{align*}
   x_k&:=(a_0,\dotsc,a_k,0,0,\dotsc) \in V
   \quad\text{and}\quad
  \\ 
   y_k&:=(a_0,\dotsc,a_{k-1},a_k+2^{-k},0,0,\dotsc) \in \WM \setminus V.
 \end{align*}
 for all $k \in \IN$. 
\begin{mycases}
 \item[``$k=0$'':]
  Set $a_0:=0$. Then $x_0=0^\omega \in V$ and $y_0=10^\omega \not\in V$.
 \item[``$k-1 \to k$'':]
  Assume that $a_0,\dotsc,a_{k-1}$ are already constructed
  with $x_{k-1} \in V \not\ni y_{k-1}$.
 \\
  Let $z:=(a_0,\dotsc,a_{k-1},1,0,0,\dotsc)$.
  As $1\leq \nm{z} \leq 2$,
  we have $z \in \WM \setminus  B_\WM(0^\omega;1)$.
  Therefore there is some $b \in M_k \setminus \{1\}$ with
  $(a_0,\dotsc,a_{k-1},b,0,0,\dotsc) \in V$ and
  $(a_0,\dotsc,a_{k-1},b+2^{-k},0,0,\dotsc) \notin V$.
  Since $(\sum_{i=0}^{k-1} a_i ) + b+2^{-k}<\infty$,
  the number $a_k:=b$ satisfies the requirements.
\end{mycases}

\noindent
 We set $x_\infty:=(a_0,a_1,\dotsc) \in \Tprod_{i \in \IN}M_i$.
 Clearly, both the sequences $(x_k)_k$ and $(y_k)_k$ converge to
 $x_\infty$ in $\IR^\IN$.
 By Lemma~\ref{l:ConvInEll} they converge to $x_\infty$ in the space $\WM$
 as well, because
 \[
   \nm{x_\infty}
   = \Tsum_{i=0}^\infty a_i
   = \lim_{m \to \infty} \nm{y_m}
   = \lim_{m \to \infty} \nm{x_m}
   \leq 1
   <\infty \,.
 \]
 As $\WM \setminus V$ is closed, we have $x_\infty \notin V$.
 Thus $V$ is not closed.
 We conclude that the ball $B_\WM(0^\omega;1)$ does not contain any
 clopen neighbourhood of $0^\omega$.
\end{Proof}

The ball $B_\WM(0^\omega;1)$ is a basic open 
of the metrisable topology of $\WM$.
Hence the complement $\WM \setminus B_\WM(0^\omega;1)$
is closed and, as $\WM$ is a metric space, even functionally closed.
Lemma~\ref{l:unitball:notZO} entails that $\WM$ is not zero-dimensional
and that $\WM \setminus B_\WM(0^\omega;1)$
is not an intersection of clopen sets.

\begin{Lemma}\label{l:WM:notNC}
 The space $\WM$ is a Polish space that is not zero-dimensional
 and that does not satisfy Normann's condition.
\end{Lemma}

In Section~\ref{sec:WM:into:INNN} we will prove that $\WM$
is a retract of the $\QCB$-space $\INNN$.
Hence $\WM$ is homeomorphic to a closed subspace of $\INNN$.
By being topologised by the sequentialisation of the compact-open topology
on the set $\INNN$, the $\QCB$-exponential $\INNN$ is
the sequential coreflection of some zero-dimensional space.
This property is inherited by closed subspaces.
We conclude that $\WM$ is the sequential coreflection of
some zero-dimensional topological space.


\section{Embedding $\WM$ into $\IN^{\IN^\IN}$ as a retract}
\label{sec:WM:into:INNN}

In this section we show that $\WM$ is a retract of
the $\QCB$-space $\INNN$.
We do this via the space $\Wzwei^{\IN \times \IF}$,
where $\IF$ denotes the countable metric \emph{fan}
and $\Wzwei$ denotes the two point discrete space with points $\w{0},\w{1}$.


\subsection{The countable fan $\IF$}
\label{sub:fanspace}

The fan space $\IF$ is the ``smallest''
non-locally-compact metrisable space,
in the sense that it embeds
into every non-locally-compact metrisable space as a closed subspace.
Our version of the countable fan has $\IN^2 \cup \{(\infty,\infty)\}$
as underlying set and its topology is induced by the unique metric $d_\IF$
that satisfies
\[
 d_\IF\big( (a,b), (\infty,\infty) \big)=2^{-a}
 \;\;\text{and}\;\;
 d_\IF\big( (a,b),(a',b') \big)=\max\big\{ 2^{-a},  2^{-a'} \big\}
\]
for $a,b,a',b'\in \IN$ with $(a,b) \neq (a',b')$.
So every point apart from $(\infty,\infty)$ is an isolated point in $\IF$.
Furthermore, $(a_n,b_n)_n$ converges to $(\infty,\infty)$
in $\IF$ if, and only if, $\lim\limits_{n \to \infty} a_n=\infty$.

By being a zero-dimensional Polish spaces,
the product $\IN \times \IF$ is a retract of the Baire space $\IN^\IN$
in $\QCB$, i.e.,
there are continuous functions $e\colon  \IF \times \IN \to \IN^\IN$
and $r\colon \IN^\IN \to \IN \times \IF$ satisfying
$r \circ e=\id_{\IN \times \IF}$
(see \cite[Section 3.3]{Sch:NNNnotT3} for an explicit construction).

 
\subsection{The zero-dimensional space $\WMProd$}

For any $i \in \IN$, 
we endow the finite set $M_i$ with its discrete topology
and denote the corresponding finite discrete space by $\WM_i$.
The product $\WMProd:=\Tprod_{i \in \IN} \WM_i$
is a zero-dimensional compact metrisable space.
Since the spaces $\WM_i$ are discrete subspaces of $\IR$,
a sequence $(x_n)_n$ of elements of the set $\Tprod_{i \in \IN} M_i$
converges in the space $\IR^\IN$ to some point $x_\infty$ 
if, and only if, it does in the space $\WMProd$.
We obtain by Lemma~\ref{l:ConvInEll}:

\begin{Lemma}\label{l:ConvInWM}
 Let $(x_n)_n$ be a sequence in $\WM$ and let $x_\infty \in \WM$.
 Then $(x_n)_n$ converges in $\WM$ to $x_\infty$
 if, and only if, (a) and (b) hold:
 \begin{itemize}
  \item[(a)]
   For all $i \in \IN$ there is some $n_i \in \IN$ with
   $x_n(i)=x_\infty(i)$ for all $n \geq n_i$.
  \item[(b)]
   The sequence $(\nm{x_n})_n$ converges to $\nm{x_\infty}$ in $\IR$.
 \end{itemize}
\end{Lemma}

Hence the injection
$\id\colon \WM \to \WMProd$ is sequentially continuous
and thus topologically continuous, as $\WM$ is metrisable.
This implies that the topology of $\WM$ is finer than the subspace topology
on the set $\setof\WM$ induced by the topology of $\WMProd$.
In fact, it is strictly finer than the subspace topology,
because the sequence $(0^{n+1}\tfrac{1}{2} 0^\omega)_n$ converges
in $\WMProd$ to $0^\omega$, but not in $\WM$.


\subsection{An embedding of $\WM$ into $\WMProd \times \Wzwei^{\IN \times \IF}$}

We now start to prove that $\WM$ is a retract of the $\QCB$-product
$\WMProd \times \Wzwei^{\IN \times \IF}$.
First we define two functions
$f\colon \WMProd \times \IN \times \IN^2 \to \{\w{0},\w{1}\}$
and
$g\colon \WM \to \Wzwei^{\IN \times \IF}$ by
\begin{align}
  &  f(x,k,a,b) := \left\{
    \begin{array}{ll}
      \w{0} & \text{if $\Tsum_{i=a}^{a+b} x(i) \leq 2^{-k}$}
      \\
      \w{1} & \text{otherwise}
    \end{array}  \right.
   \\ 
  &  g(y)(k,a,b) := f(y,k,a,b)
   \;\;\text{and}\;\;
   g(y)(k,\infty,\infty):= \mathtt{0}
\end{align}
for all $x \in \WMProd$, $y \in \WM$ and $k,a,b \in \IN$.

\pagebreak[3]
\begin{Lemma}\label{l:f:g:continuous}
 For any $y \in \WM$, the function $g(y)\colon \IN \times \IF \to \Wzwei$
 is continuous.
 Moreover, $f$ and $g$ are continuous.
\end{Lemma}

\begin{Proof}
\begin{enumerate}
 \item
  Let $(x_n,k_n,a_n,b_n)_n$ converge
  to $(x_\infty,k_\infty,a_\infty,b_\infty)$
  in $\WMProd \times \IN \times \IN^2$.
  Then there is some $n_0 \in \IN$ such that
  \[
   x_n(i)=x_{n_0}(i)
   ,\quad
   k_n=k_\infty
   \quad\text{and}\quad
   (a_n,b_n)=(a_\infty,b_\infty)
  \]
  for all $n \geq n_0$ and
  all $i \in \{a_\infty,\dots, a_\infty+b_\infty\}$.
  By the definition of $f$, we have
  $f(x_n,k_n,a_n,b_n)=f(x_\infty,k_\infty,a_\infty,b_\infty)$
  for all $n \geq n_0$.
  So $f$ is continuous.
 \item
  By the cartesian closedness of $\QCB$ it suffices to show
  the continuity of the function
  $G\colon \WM \times \IN \times \IF \to \Wzwei$
  defined by
  \[
    G(x,k,a,b)= f(x,k,a,b)
    \;\;\text{and}\;\;
    G(x,k,\infty,\infty):=\mathtt{0} \;.
  \]
  Let $(x_n,k_n,a_n,b_n)_n$ converge to
  $(x_\infty,k_\infty,a_\infty,b_\infty)$ in
  $\WM \times \IN \times \IF$.
  If $(a_\infty,b_\infty) \in \IN^2$, then
  the sequence $\big(G(x_n,k_n,a_n,b_n)\big){}_n$
  converges to $G(x_\infty,k_\infty,a_\infty,b_\infty)$
  by the continuity of $f$ and by the fact that the topology of $\WM$
  is finer than the subspace topology on $\WM$
  induced by the topology of $\WMProd$.
  \\
  Now let $(a_\infty,b_\infty)=(\infty,\infty)$.
  Since $\Tsum_{i=0}^\infty x_\infty(i) < \infty$, there is some
  $m \in \IN$ with $\Tsum_{i=m}^\infty x_\infty(i)< 2^{-k_\infty-1}$.
  There exists some $n_0 \in \IN$ such that,
  for all $n \geq n_0$ and all $i < m$,
  \[
   x_n(i)=x_\infty(i)
   ,\;\;
   \big| \nm{x_n} - \nm{x_\infty} \big| < 2^{-k_\infty-1}
   ,\;\;
   k_n=k_\infty
   \;\;\text{and}\;\;
   a_n \geq m
   \,.
  \]
  Then all $n \geq n_0$ with $(a_n,b_n)\neq (\infty,\infty)$
  satisfy
  \begin{align*}
    \Tsum_{i=a_n}^{a_n+b_n} x_n(i)
    &\leq \Tsum_{i=m}^{\infty} x_n(i)
     = \nm{x_n} - \Tsum_{i=0}^{m-1} x_n(i)
     = \nm{x_n} - \Tsum_{i=0}^{m-1} x_\infty(i)
    \\ 
     &= \nm{x_n} - \nm{x_\infty} + \Tsum_{i=m}^\infty x_\infty(i)
      < 2 \cdot 2^{-k_\infty-1}
      = 2^{-k_n} \,.
  \end{align*}
  This implies
  $G(x_n,k_n,a_n,b_n)=\w{0} =G(x_\infty,k_\infty,a_\infty,b_\infty)$
  for all $n \geq n_0$.
  Hence $G$ is sequentially continuous
  and thus topologically continuous as a function
  from the $\QCB$-product $\WM \times \IN \times \IF$
  to the two-point discrete space $\Wzwei$.
  We conclude that $g$ is a continuous function
  into the space $\Wzwei^{\IN \times \IF}$.
 \end{enumerate}  
\end{Proof}

Using the continuous map $g$, we define the function
$e_\WM\colon \WM \to \WMProd \times \Wzwei^{\IN \times \IF}$ by
\begin{equation}\label{eq:Definition:eWM}
  e_\WM(x):=\big(x,g(x)  \big) 
\end{equation}
for all $x \in \WM$.
From Lemmas~\ref{l:ConvInWM} and~\ref{l:f:g:continuous} we obtain:

\begin{Lemma}
 The function $e_\WM$ is injective and continuous.
\end{Lemma}


\subsection{Construction of the retraction map from
   $\WMProd \times \Wzwei^{\IN \times \IF}$ to $\WM$}
  
To construct a retraction map pertaining to the section $e_\WM$,
we define for any $m \in \IN$ the set $C_m$ by
\begin{align*}
 C_m:=
 \Big\{ (x,h) \in \WMProd \times \Wzwei^{\IN \times \IF} \;\Big|\;
        &  h(k,\infty,\infty)=\w{0} \text{ and}
    \\	&\quad  h(k,a,b)= f(x,k,a,b)
	 \;\;\text{for all $k,a,b \leq m$}
 \Big\}
 \,.
\end{align*}
Using the decreasing sequence $(C_m)_m$, we construct a function
$r_\WM$ from $\WMProd \times \Wzwei^{\IN \times \IF}$ to
$\Tprod_{i \in \IN} M_i$ by
\[
 r_\WM(x,h)(m):=\left\{
 \begin{array}{cl}
  x(m) & \text{if $(x,h) \in C_m$}
  \\
   0   & \text{otherwise}
 \end{array}\right.
\]
for all $(x,h) \in \WMProd \times \Wzwei^{\IN \times \IF}$
and $m \in \IN$.
It follows immediately from the definitions
that the image of $e_\WM$ lies in $\Tbigcap_{m \in \IN} C_m$
and that $r_\WM(e_\WM(x))=x$ holds for every $x \in \WM$.

We show that $r_\WM$ maps into the space $\WM$.
This fact allows us to consider $r_\WM$ henceforth as a function
of the form $\WMProd \times \Wzwei^{\IN \times \IF} \to \WM$.

\pagebreak[3]
\begin{Lemma}\label{l:rWM:properties}\quad\nopagebreak
\begin{enumerate}\vspace*{-0.5ex}
 \item\label{en:rWM:h:property}
  Let $(x,h) \in \WMProd \times \Wzwei^{\IN \times \IF}$.
  Let $(a,k) \in \IN^2$ such that $a \geq k$ and
  $h(k,a,b)=\w{0}$ for all $b \in \IN$.
  Then we have $\Tsum_{i=a}^\infty r_\WM(x,h)(i) \leq 2^{-k}$.
 \item\label{en:bigcap:C_m:image:eM}
  The image of $e_\WM$ is equal to 
  the intersection $\Tbigcap_{m \in \IN} C_m$.
 \item\label{en:rWM.eWM=id}
  For every $(x,h) \in \WMProd \times \Wzwei^{\IN \times \IF}$,
  we have $r_\WM(x,h) \in \WM$.
 \item\label{en:Cm:clopen} 
  For every $m \in \IN$, the set $C_m$
  is clopen in $\WMProd \times \Wzwei^{\IN \times \IF}$.  
\end{enumerate}
\end{Lemma}

\pagebreak[3]
\begin{Proof}
\begin{enumerate}\vspace*{-0.5ex}
 \item 
  By induction on $m$ we show
  $\Tsum_{i=a}^m r_\WM(x,h)(i) \leq 2^{-k}$ for all $m \geq a$.
  \begin{mycases}
   \item[``$m=a$'':]
    If $(x,h) \in C_a$, then $f(x,k,a,0)=h(k,a,0)=\w{0}$
    and $r_\WM(x,h)(a)=x(a)=\Tsum_{i=a}^a x(i) \leq 2^{-k}$.
    Otherwise we have $r_\WM(x,h)(a)=0 \leq 2^{-k}$.
   \item[``$m>a$'':] 
    If $(x,h) \in C_m$, then we have $r_\WM(x,h)(i)=x(i)$ for all $i \leq m$
    and $f(x,k,a,m-a)=h(k,a,m-a)=\w{0}$, because $k \leq m$.
    This implies
    $\Tsum_{i=a}^m r_\WM(x,h)(i)=\Tsum_{i=a}^m x(i) \leq 2^{-k}$.
    \\
    Otherwise $r_\WM(x,h)(m)$ is equal to $0$ and
    the induction hypothesis yields
    $\Tsum_{i=a}^m r_\WM(x,h)(i)
       = \Tsum_{i=a}^{m-1} r_\WM(x,h)(i) \leq 2^{-k}$.
  \end{mycases}
  We conclude $\Tsum_{i=a}^\infty r_\WM(x,h)(i) \leq 2^{-k}$.
 \item
  We have already noticed $e_\WM[\WM] \subseteq \Tbigcap_{m \in \IN} C_m$.
  To show ``$\supseteq$'',
  let $(x,h) \in \Tbigcap_{m \in \IN} C_m$.
  \\
  The continuity of $h$
  implies that there is some $a_0 \in \IN$ such that
  $h(0,a,b)=h(0,\infty,\infty)=\w{0}$
  for all $a \geq a_0$ and $b \in \IN$.
  From \eqref{en:rWM:h:property} we obtain
  \begin{align*}
   \nm{x}
   =
   \Tsum_{i=0}^ \infty x(i)
   =
   \Tsum_{i=0}^\infty r_\WM(x,h)(i)
   &=
   \Tsum_{i=a_0}^\infty r_\WM(x,h)(i) + \Tsum_{i=0}^{a_0-1} r_\WM(x,h)(i) 
   \\
   &\leq 2^0 +\Tsum_{i=0}^{a_0-1} r_\WM(x,h)(i)
   <\infty,
  \end{align*}
  hence $x \in \WM$.
  This allows us to apply $g$ and $e_\WM$ to $x$.
  Since $(x,h) \in \Tbigcap_{m\in\IN} C_m$, we have $g(x)=h$ and $e_\WM(x)=(x,h)$.
  Therefore $(x,h)$ lies in the image of $e_\WM$.
 \item
  Clearly, $r_\WM(x,h)(i) \in M_i$ holds for every $i \in \WM$.
  If $(x,h) \in e_\WM[\WM]$, then we have
  $(x,h) \in \Tbigcap_{m \in \IN} C_m$ and thus $r_\WM(x,h)=x \in \WM$.
  Otherwise, if $(x,h) \notin e_\WM[\WM]$,
  then by~\eqref{en:bigcap:C_m:image:eM}
  there is some $m \in \IN$ with $(x,h) \notin C_m$.
  This implies $r_\WM(x,h)(i)=0$ for all $i \geq m$.
  Hence $\nm{r_\WM(x,h)} < \infty$ and thus $r_\WM(x,h) \in \WM$.
 \item
  Let $k,a,b \in \{0,\dotsc,m\}$.
  For every $s \in \{\w{0},\w{1}\}$ the set
  \[
      D_s:= \big\{ x \in \WMProd \,\big|\, f(x,k,a,b)=s \big\}
  \]
  is clopen in $\WMProd$ by the continuity of $f$
  (see Lemma~\ref{l:f:g:continuous}).
  Moreover, the set
  \[
   E_s:= \big\{ h \in \Wzwei^{\IN \times \IF} \,\big|\,
          h(k,\infty,\infty)=\w{0} \;\;\text{and}\;\; h(k,a,b)=s \big\}
  \]
  is clopen w.r.t.\ the compact-open topology
  and thus w.r.t.\ the sequential topology on
  $\Wzwei^{\IN \times \IF}$, because the latter is finer than the former.
  Hence the set
  \[
   \big\{ (x,h) \in \WMProd \times \Wzwei^{\IN \times \IF} \;\big|\;
       h(k,\infty,\infty)=\w{0} \;\;\text{and}\;\; h(k,a,b)= f(x,k,a,b)
   \big\}	
  \]
  is clopen in $\WMProd \times \Wzwei^{\IN \times \IF}$
  by being equal to
  $(D_{\w{0}} \times E_{\w{0}}) \cup (D_{\w{1}} \times E_{\w{1}})$.
  Therefore $C_m$ is clopen by being a finite intersection of clopen sets.
\end{enumerate}
\end{Proof}

We need the following lemma about converging sequences
in the $\QCB$-exponential $\Wzwei^{\IN \times \IF}$.
It can be easily deduced from the fact that the convergence relation
of $\QCB$-exponentials is continuous
convergence\footnote{A sequence $(f_n)_n$ of continuous functions
  between two sequential spaces $\WX$ and $\WY$
  is said to \emph{converge continuously}
  to a continuous function $f_\infty\colon \WX \to \WY$,
  if $(f_n(x_n))_n$ converges to $f_\infty(x_\infty)$ in $\WY$,
  whenever $(x_n)_n$ converges to $x_\infty$ in $\WX$.},
  see \cite{ELS:CCGS,Sch:phd}.

\pagebreak[3]
\begin{Lemma}\label{l:to:IF:convergence:prop} 
  Let $(h_n)_n$ converge to $h_\infty$ in $\Wzwei^{\IN \times \IF}$.
  Then for every $k \in \IN$ there exists
  some $m \in \IN$
  with
  $ h_n(k,a,b)= h_\infty(k,\infty,\infty) $
  for all $n \geq m$ (including $n=\infty$),
  all $a \geq m$ and all $b \in \IN$.
\end{Lemma}

Now we are able to show that $r_\WM$ is a retraction map.

\pagebreak[3]
\begin{Proposition}\label{p:WM:is:retract}
  The space $\WM$ is a retract of
  $\WMProd \times \Wzwei^{\IN \times \IF}$ in $\QCB$.
  The functions $e_\WM$ and $r_\WM$ form
  a section-retraction-pair.
\end{Proposition}

\begin{Proof}
  We have already observed $r_\WM \circ e_\WM=\id_\WM$
  and the continuity of $e_\WM$.
  \\
  It remains to prove the continuity of $r_\WM$.
  Let $(x_n,h_n)_n$ converge to $(x_\infty,h_\infty)$.
  We set $z_n:=r_\WM(x_n,h_n)$ for all $n \in \IN \cup\{\infty\}$.
  For every $m \in \IN$, $(z_n(m))_n$ converges to $z_\infty(m)$ in $\WM_m$,
  because $C_m$ is clopen by Lemma~\ref{l:rWM:properties}
  and $(x_n(m))_n$ converges to $x_\infty(m)$ in $\WM_m$.
  We consider two cases.
  \begin{itemize}
   \item[(a)] Let $(x_\infty,h_\infty) \in e_\WM[\WM]$.
    Then $z_\infty=x_\infty$, 
    because $(x_\infty,h_\infty) \in \Tbigcap_{m \in \IN} C_m$.
    \\
    To show $\lim\limits_{n \to \infty} \nm{z_n - z_\infty}=0$,
    let $k \in \IN$.
    Then $h_\infty(k,\infty,\infty)=\w{0}$,
    because $(x_\infty,h_\infty) \in C_k$.
    By Lemma~\ref{l:to:IF:convergence:prop}
    there exists some $m \in \IN$ such that
    \[
     h_n(k,a,b)=h_\infty(k,\infty,\infty)=\w{0}
    \]
    for all $n \geq m$ (including $n=\infty$), $a \geq m$ and $b \in \IN$.
    We set $c:=\max\{m,k\}$.
    Since $(z_n)_n$ converges pointwise to $z_\infty$,
    there is some $n_1 \in \IN$ with $z_n(i)=z_\infty(i)$
    for all $n \geq n_1$ and $i <c$.
    By Lemma~\ref{l:rWM:properties}\eqref{en:rWM:h:property},
    all $n \geq \max\{m,n_1\}$ satisfy
    \begin{align*}
     \nm{z_n - z_\infty}
     &=\big| \Tsum_{i=0}^\infty z_n(i)
             - \Tsum_{i=0}^\infty z_\infty(i) \big|
      =\big| \Tsum_{i=c}^\infty z_n(i)
             - \Tsum_{i=c}^\infty z_\infty(i) \big|
     \\
     &\leq \max\big\{ \Tsum_{i=c}^\infty z_n(i),\;
                 \Tsum_{i=c}^\infty z_\infty(i) \big\}
      \leq 2^{-k}.
    \end{align*}
    We conclude that $(z_n)_n$ converges to $z_\infty$ in $\leins$
    and thus in the subspace $\WM$.
   \item[(b)]
    Let $(x_\infty,h_\infty) \notin e_\WM[\WM]$.
    By Lemma~\ref{l:rWM:properties}\eqref{en:bigcap:C_m:image:eM}
    there is some $m \in \IN$ with $(x_\infty,h_\infty) \notin C_m$.
    Since $(z_n)_n$ converges pointwise to $z_\infty$
    and $C_m$ is closed, there is some $n_0 \in \IN$ such that,
    for all $n \geq n_0$ and $i<m$,
    \[
     (x_n,h_n) \notin C_m
     \quad\text{and}\quad
     z_n(i)=z_\infty(i)
     \,.
    \]
    For all $n \geq n_0$ we have $z_n=z_\infty$,
    because $z_n(i)=0=z_\infty(i)$ holds for all $i \geq m$.
    Therefore $(z_n)_n$ converges to $z_\infty$ in $\WM$.
  \end{itemize}
\end{Proof}


\subsection{Establishing $\WM$ as a retract of $\INNN$}

As $\WMProd$ is a  zero-dimensional compact metrisable space
without isolated points,
$\WMProd$ is homeomorphic to the Cantor space $\Wzwei^\IN$
by Theorem~7.4 in \cite{Kech:descriptive}.
Moreover, the product $2^\IN \times \Wzwei^{\IN \times \IF}$
is homeomorphic to $\Wzwei^{\IN \times \IF}$,
hence $\WMProd \times \Wzwei^{\IN \times \IF}$
is homeomorphic to $\Wzwei^{\IN \times \IF}$.
Since $\IN \times \IF$ is a retract of $\IN^\IN$
(see Section~\ref{sub:fanspace}),
$\Wzwei^{\IN \times \IF}$ is a retract of $\IN^{\IN^\IN}$.
We obtain by Proposition~\ref{p:WM:is:retract}:

\begin{Proposition}\label{p:WM:retract:2toINxIF:NNN}
 The Polish space $\WM$ is a retract of $\Wzwei^{\IN \times \IF}$
 and of $\IN^{\IN^\IN}$.
\end{Proposition}


\section{The main result}
\label{sec:NNN:notNC}

To establish our main result that $\INNN$ does not satisfy
Normann's condition, it remains to verify that forming
retracts preserves Normann's condition.

\begin{Lemma}\label{l:retracts:inherit:NP}
 Let $\WX$ be a retract of some qcb-space $\WY$.
 If $\WY$ satisfies Normann's condition,
 then so does $\WX$.
\end{Lemma}

\begin{Proof}
 Let $e\colon \WX \to \WY$ and $r\colon \WY \to \WX$
 be continuous functions with $r \circ e = \id_\WX$.
 Let $A$ be a functionally closed subset of $\WX$.
 Then $r^{-1}[A]$ is a functionally closed subset of $\WY$.
 As $\WY$ satisfies Normann's condition,
 there is a family of clopen subsets $(C_i)_{i \in I}$ of $\WY$
 with $\Tbigcap_{i \in I} C_i= r^{-1}[A]$.
 By the continuity of $e$, the sets $e^{-1}[C_i]$ are clopen in $\WX$.
 One easily verifies $A=\Tbigcap_{i \in \IN} e^{-1}[C_i]$.
 Therefore $A$ is an intersection of clopen sets of $\WX$.
\end{Proof}

Hence Lemma~\ref{l:WM:notNC} and
Proposition~\ref{p:WM:retract:2toINxIF:NNN} imply our main result
stating that $\IN^{\IN^\IN}$ contains functionally closed sets
that are not intersections of clopens.

\begin{Theorem}\label{th:INNN:not:NC}
 The space $\IN^{\IN^\IN}$ does not satisfy Normann's condition.
\end{Theorem}

Analogously, neither $\WMProd \times \Wzwei^{\IN \times \IF}$
nor $\Wzwei^{\IN \times \IF}$ satisfies Normann's condition.
By induction on $k$ one can show that $\Nk{k}$ is a retract
of $\Nk{k+1}$. We conclude by Lemma~\ref{l:retracts:inherit:NP}:

\begin{Corollary}\label{c:Nkk:notNC}
 For every $k \geq 2$, the sequential space $\Nk{k}$ of Kleene-Kreisel functionals
 of level $k$ does not satisfy Normann's condition.
\end{Corollary}


\section{An application in Computable Analysis}
\label{sec:application:CompAna}

In \cite{BES:paradigms}, Bauer, Escard{\'o} and Simpson formalised
two approaches to higher type computation over the reals numbers
in functional programming
by defining two ``real'' objects in the category $\Equ$ of equilogical spaces
\cite{BBS:Equi}.
The first object, $\RE$, models the \emph{external reals} describing
the approach of introducing the reals as an own datatype.
The object $\RI$ models the concept of representing real numbers
via infinite streams.
These reals are called \emph{internal reals}.

Using exponentiation in the cartesian closed category $\Equ$,
the application of the natural recursion formulae 
\begin{alignat}{3}
 \REk{0}&:= \RE  &&\;\;\text{and}\;\;& \REk{k+1}&:=\RE^{\REk{k}} \\
 \RIk{0}&:= \RI  &&\;\;\text{and}\;\;& \RIk{k+1}&:=\RI^{\RIk{k}}
\end{alignat}
yields two hierarchies of functionals over the real numbers.
The hierarchy of the underlying sets of the sequence $(\REk{k})_k$
is called the \emph{extensional hierarchy}.
The underlying sets of $(\RIk{k})_k$
form the \emph{intensional hierarchy}.

The natural question arises whether the two hierarchies of functionals
coincide. This question is referred to as the \emph{Coincidence Problem}.
From \cite{BES:paradigms} we know that both hierarchies agree up to level $2$.
Normann's equivalence result
(Theorem 4.17 and 5.5 in \cite{Nor:comparing})
states that the two hierarchies agree on level $k+1$ if, and only if,
every functionally closed subset of the Kleene-Kreisel space $\Nk{k}$
is an intersection of clopen sets.
Therefore our main result (Theorem~\ref{c:Nkk:notNC})
along with Corollary~\ref{c:Nkk:notNC}
solves the Coincidence Problem negatively.

\begin{Theorem}
 The extensional hierarchy and the intensional hierarchy
 of functionals  over the reals
 do not coincide from level $3$ on.
\end{Theorem}
Hence both hierarchies disagree from the first previously unknown level on.
It is known that the extensional hierarchy coincide
with the sequential hierarchy. The latter is formed by the underlying
sets of the sequence
\[
 \RSk{0}:=\IR \;\;\text{and}\;\;\RSk{k+1}:=\IR^{\RSk{k}}
\]
formed in the category $\QCB$ (equivalently in $\Seq$ or $\kHaus$).
So Theorem 5.1 states there is a continuous functional
$F\colon \IR^{\IR^\IR} \to \IR$ that is not an element of 
the space $\RIk{3}$.

 
\pagebreak[3]

\bigskip
\noindent
\texttt{Fakult\"at f\"ur Informatik, Universit\"at der Bundeswehr, Munich, Germany}
\\
\texttt{\emph{Email:}\;matthias.schroeder@unibw.de}


\end{document}